\documentclass[twoside,onecolumn]{article}
\usepackage{amsmath,amssymb,epsfig}
\title{\textsc{On the propagation of an optical wave in a photorefractive medium}}

\author{\small B. BID\'EGARAY-FESQUET\\
\it\small Laboratoire Jean Kuntzmann, \\
\it\small CNRS UMR 5224 et Grenoble Universit\'es, \\ 
\it\small B.P. 53, 38041 Grenoble Cedex 9, France\\
\it\small Brigitte.Bidegaray@imag.fr\\[5mm]
\small J.-C. SAUT\\
\it\small UMR de Math\'ematiques, Universit\'e Paris-Sud,\\ 
\it\small 91405 Orsay, France\\
\it\small Jean-Claude.Saut@math.u-psud.fr}

\date{}

\begin{document}

\maketitle

\begin{abstract}
The aim of this paper is first to review the derivation of a model describing
the propagation of an optical wave in a photorefractive medium and to present
various mathematical results on this model: Cauchy problem, solitary waves.
\end{abstract}



\renewcommand{\thefootnote}{\fnsymbol{footnote}}

\footnote[0]{Electronic version of an article published as \\	
Mathematical Models and Methods in Applied Sciences (M3AS) \\
Volume 17, Issue 11, 2007, 1883--1904, DOI 10.1142/S0218202507002509 \\
\copyright World Scientific Publishing Company \\
{\texttt{http://www.worldscinet.com/cgi-bin/details.cgi?id=voliss:m3as{\textunderscore}1711\&type=toc}}}

%%%%%%%%%%%%%% AUTHORS' MACROS %%%%%%%%%%%%%% 

\def\ds{\displaystyle}
\def\ph{\varphi}
\def\eps{\varepsilon}
\def\bfE{\mathbf{E}}
\def\bfJ{\mathbf{J}}
\def\bfc{\mathbf{c}}
\def\bfe{\mathbf{e}}
\def\bfk{\mathbf{k}}
\def\bfr{\mathbf{r}}
\def\bfx{\mathbf{x}} 
\def\bbR{\mathbb{R}}
\def\calC{\mathcal{C}}
\def\calD{\mathcal{D}}
\def\calF{\mathcal{F}}
\def\calS{\mathcal{S}}
\renewcommand{\div}{\operatorname{div}}
\newcommand{\tens}[1]{\widehat{#1}}

\newtheorem{theorem}{Theorem}
\newtheorem{lemma}[theorem]{Lemma}
\newtheorem{corollary}[theorem]{Corollary}
\newtheorem{proposition}[theorem]{Proposition}

\def\theequation{\arabic{section}.\arabic{equation}}
\newcommand{\refcite}[1]{\ \cite{#1}}
\def\proof{\noindent {\bf Proof.\ }}
\def\endproof{$\Box$}

\def\remark{\noindent\textbf{Remark.}\ }

\newlength{\fntxvi} \newlength{\fntxvii}
\newcommand{\chemical}[1]
{{\fontencoding{OMS}\fontfamily{cmsy}\selectfont
  \fntxvi\the\fontdimen16\font
  \fntxvii\the\fontdimen17\font
  \fontdimen16\font=3pt \fontdimen17\font=3pt
  $\mathrm{#1}$
  \fontencoding{OMS}\fontfamily{cmsy}\selectfont
  \fontdimen16\font=\fntxvi \fontdimen17\font=\fntxvii }}

%%%%%%%%%%%%%% END AUTHORS' MACROS %%%%%%%%%%%%%% 

\maketitle

\section{Introduction}

A modification of the refraction index in \chemical{LiNbO_3} or
\chemical{LiTaO_3} crystals has been observed in the 1960s and first 
considered as a drawback. This photo-induced varia\-tion of the index is called 
the photorefractive effect and occurs in any electro-optical or photoconductive 
crystal. Applications have been found in the 1970s--1980s to real-time 
signal processing, phase conjugation, or amplification of beams or images. 

In this paper we are interested in deriving a not-too-simple but tractable 
math\-ematical model for the propagation of light in such materials. Solitonic 
propagation is one of our concern but we focus here on initial value 
problems. A very complete review of solitonic propagation in photorefractive 
media may be found in \refcite{delre-crosignani-diporto}. Our derivation 
follows the same guidelines as theirs but point out the different 
approximations made for future mathe\-matical studies. 

The outline of the paper is the following. In Sec. \ref{Sec_Model} we 
first derive the Kukhtarev model for the material and then couple it to a wave 
propagation model for light to obtain a complete set of equations. A 1D model
is obtained keeping only one of the two transverse space variables. This is a 
saturated nonlinear Schr\"odinger equation, the mathematical theory of which 
is addressed in Sec. \ref{Sec_saturated} in arbitrary dimension. Section 
\ref{Sec_Cauchy} is devoted to the study of the full 2D model with emphasis on
the Cauchy problem and the solitary wave solutions. 

\section{Derivation of the Model}
\label{Sec_Model}

\subsection{The photorefractive effect}

The propagation of an optical wave in insulating or semi-insulating 
electro-optical crystals induces a charge transfer. The new distribution of 
charges induces in turn an electric field which produces a variation of the 
refraction index. The main characteristics of this effect are the following: 
(1) Sensibility to energy (and not to the electric field), 
(2) Nonlocal effect (charge distributions and the electric field are not 
located at the same position),
(3) Inertia (charges need a certain time to move), 
(4) Memory and reversibility (in the dark the space charge, and therefore
the index variation, is persistent but an uniform light redistributes 
uniformly all charges --- this yields applications to holography). 

The sensibility to energy reminds us of Kerr media yielding the classical 
cubic nonlinear Schr\"odinger (NLS) equation. The nonlocal effects will of 
course complicate the mathematical analysis compared to NLS equations, but the 
general ideas will be the same. In our final model, inertial effects will be 
neglected since time is removed from the material equations. Memory and 
reversibility effects involve ion displacement in materials like 
\chemical{Bi_2TeO_5}, which we will not take into account in the present 
study. 

\subsection{The Kukhtarev model}

The physical modeling of the photorefractive effect assumes that charges are 
trapped in impurities or defaults of the crystal mesh. We chose here to derive 
the model only in the case when charges are electrons. Some materials like 
semi-conductors necessitate to model both electrons and holes. Therefore we 
restrict our study to insulating media. 

\subsubsection{Charge equation}

Electrons come from donor sites with density $N_{\rm D}$. This density is 
supposed to be much greater than that of the acceptor sites (impurities) which 
we denote by $N_{\rm A}$. The density of donor sites which are indeed ionized 
is $N_{\rm D}^+$ and we of course have 
$N_{\rm D}^+\leq N_{\rm A} \ll N_{\rm D}$. Local neutrality, i.e. no electrons 
in the conduction band, corresponds to the relation $N_{\rm D}^+=N_{\rm A}$. 
The total charge is given by
\begin{equation}
\label{Eq_Charge} %2.1
\rho = e(N_{\rm D}^+-N_{\rm A}-n_{\rm e}),
\end{equation}
where $e$ is the electron charge and $n_{\rm e}$ the electron density.

\subsubsection{Evolution of ionized donor sites}

Photoionization and recombination affect the density of ionized donor sites.
Photoionization is proportional to the density of not ionized donor sites 
($N_{\rm D}-N_{\rm D}^+$). In the dark it is proportional to a thermal 
excitation rate $\beta$ but is also sensitive to light intensity $I_{\rm em}$ 
with a photoexcitation coefficient $s$. Recombination is proportional to the 
density of electrons and occurs over a time scale 
$\tau=1/(\gamma_{\rm r}N_{\rm D}^+)$ which does not depend on $n_{\rm e}$ if 
the excitation rate is low, therefore the total evolution of ionized donor 
sites is
\begin{equation}
\label{Eq_Donor} %2.2
\partial_t N_{\rm D}^+ 
= (\beta+sI_{\rm em})(N_{\rm D}-N_{\rm D}^+) 
- \gamma_{\rm r} n_{\rm e} N_{\rm D}^+.
\end{equation}

\subsubsection{Charge transport}

Now the main point is to describe the three phenomena which contribute to the 
charge transport or current density. The first phenomenon is isotropic and is 
due to thermal diffusion. It is proportional to the gradient of the electron 
density. The electron mobility is denoted by $\mu$, $T$ is the temperature and 
$k_{\rm B}$ the Boltzmann constant. The second phenomenon is drift and is 
collinear to the electric field $\bfE_{\rm tot}$. Finallly, the photovoltaic 
effect is collinear to the optical axis $\bfc$ and proportional to the 
non-ionized donor density and the field intensity with a photovoltaic 
coefficient $\beta_{\rm ph}$. The total current density is therefore
\begin{equation}
\label{Eq_Current} %2.3
\bfJ 
= e\mu n_{\rm e}\bfE_{\rm tot} 
+ \mu k_{\rm B} T\nabla n_{\rm e} 
+ \beta_{\rm ph}(N_{\rm D}-N_{\rm D}^+)\bfc I_{\rm em}.
\end{equation}

\subsubsection{Closure of the model}

The closure of the model is first based on charge conservation and the Poisson
equation:
\begin{equation}
\label{Eq_Conservation} %2.4
\partial_t \rho + \nabla\cdot\bfJ = 0, 
\end{equation}
\begin{equation}
\label{Eq_Poisson} %2.5
\nabla\cdot(\eps_0\tens\eps\bfE_{\rm sc}) = \rho.
\end{equation}
The crystal is anisotropic and this is accounted for in the relative 
permittivity~$\tens\eps$ which is a tensor. A careful analysis of the 
different electric fields has to be done. In  Poisson equation 
\eqref{Eq_Poisson}, $\bfE_{\rm sc}$ is the space charge field which is induced 
by the charge density. The total field $\bfE_{\rm tot}$ only occurs in the 
equations through its gradient (Eqs. \eqref{Eq_Conservation} and 
\eqref{Eq_Current}). Two fields are constant and disappear in the final 
equations: the photovoltaic field 
$\bfE_{\rm ph}=\beta_{\rm ph}\gamma_{\rm r}N_{\rm A}\bfc/e\mu s=E_{\rm ph}\bfc$,
and an external field $\bfE_{\rm ext}$ which is often applied in one transverse 
direction on the faces of the crystal. A last contribution to the total field 
is $\bfE$, connected to the light which propagates in the crystal and its 
description is given in Sec. \ref{Subsec_NLS}.

The set of five equations \eqref{Eq_Charge}--\eqref{Eq_Poisson} is called the 
Kukhtarev model and was first given in \refcite{kukhtarev-markow-odoluv-soskin-vinetskii}.

\subsection{Propagation of the light wave in the crystal}
\label{Subsec_NLS}

We have already introduced the relative permittivity tensor $\tens\eps$ which
plays a r\^ole in the description of the propagation of a light wave in the 
crystal via the wave equation:
\begin{equation*}
\partial_t^2 (\tens\eps \bfE) - c^2 \nabla^2\bfE = 0.
\end{equation*}
In a non-centrosymmetric crystal the preponderant nonlinear effect is the Pockels 
effect which yields the following $\bfE$-dependence for the permittivity 
tensor:
\begin{equation*}
\tens{\eps(\bfE)} 
= \tens{\eps(0)} - \tens\eps \tens{(\tens\bfr\cdot\bfE)} \tens\eps 
= n^2 - \tens\eps \tens{(\tens\bfr\cdot\bfE)} \tens\eps, 
\end{equation*} 
where $\tens\bfr$ is the linear electro-optic tensor and $n$ the mean 
refraction index. We now suppose that $\bfE$ is a space perturbation of a 
plane wave (paraxial approximation) of frequency $\omega$, wave vector $\bfk$ 
and polarization $\bfe$:
\begin{equation*}
\bfE(t,\bfx) = A(\bfx) \exp(i(\omega t-\bfk\cdot\bfx)) \bfe.
\end{equation*}
Such a wave with polarization $\bfe$ only "sees" a part of tensor 
$\tens{\eps(\bfE)}$, or equivalently a variation $\delta n$ of the refraction 
index $n$:
\begin{equation*}
\delta n 
= -\frac1{2n} [\bfe \tens\eps\ \tens{\tens{\bfr}\bfe}\ \tens\eps \bfe^*] \bfE.
\end{equation*}
Now we can write an equation for the amplitude $A$ which takes into account 
the dispersion relation $c^2|\bfk|^2 = n^2\omega^2$ and the slowly varying 
envelope approximation in the $\bfk$ direction. We denote by $\nabla_\perp$ 
the gradient in the perpendicular directions to $\bfk$ and
\begin{equation}
\label{Eq_Envelope} %2.6
\left[\nabla_\perp^2 - 2i\bfk\cdot\nabla + 2 |\bfk|^2 \frac{\delta n}{n}\right]
A(\bfx) \bfe = 0.
\end{equation}
Of course, we can consider the superposition of such waves to describe for 
example pump and probe experiments.

The system is now closed but it is impossible to solve Eqs. 
\eqref{Eq_Charge}--\eqref{Eq_Envelope}. We have to simplify them taking into 
account characteristic scales. Our description follows (or more precisely 
makes explicit the assumptions in \refcite{zozulya-anderson}) and is purely 
formal. The rigorous justification is certainly difficult and should include 
the approximations made in Sec. \ref{Subsec_NLS}.

\subsection{Characteristic values}

We first want to define a characteristic electron intensity $n_0$ by 
considering uniform solutions in space and time. Equations \eqref{Eq_Charge}, 
\eqref{Eq_Donor} and \eqref{Eq_Poisson} yield
\begin{equation*}
\rho = e(N_{\rm D}^+-N_{\rm A}-n_{\rm e}) = 0 \textrm{ and } 
0 = (\beta+sI)(N_{\rm D}-N_{\rm D}^+) - \gamma_{\rm r} n_{\rm e} N_{\rm D}^+.
\end{equation*}
With a characteristic intensity $I_0$, neglecting $\beta$ and assuming 
$n_{\rm e}\ll N_{\rm A}$, we have 
$n_0 = sI_0 (N_{\rm D}-N_{\rm A})/\gamma_{\rm r} N_{\rm A}$.

There are three characteristic times: 
(1) the characteristic lifetime of an electron (in the dark) 
$\tau_{\rm e}=1/\gamma_{\rm r}N_{\rm A}$, 
(2) the characteristic evolution time of ionized donors 
$\tau_{\rm d}=1/\gamma_{\rm r}n_0$ (and a consequence of $n_0\ll N_{\rm A}$ is 
$\tau_{\rm e}\ll\tau_{\rm d}$),
(3) the characteristic relaxation time of the electric field 
$t_0=\eps_0\eps_{\rm c}/e\mu n_0$, where $\eps_{\rm c}$ is the characteristic 
value of $\tens\eps$ along the $\bfc$ direction. It is obtained combining 
Eqs. \eqref{Eq_Current} and \eqref{Eq_Poisson} assuming there is only 
drift. If a timescale has to be kept, it is $t_0$, but we do not detail 
this point since we neglect time-dependence in the final equations. 

The Debye length $L_{\rm D}$ is the characteristic value of the field space 
variation. It is determined together with the characteristic field $E_0$. The 
Poisson equation \eqref{Eq_Poisson} yields $L_{\rm D}=\eps_0\eps_{\rm c}E_0/eN_{\rm A}$. If 
drift and isotropic diffusion have the same order, $E_0=k_{\rm B}T/eL_{\rm D}$ and 
therefore
\begin{equation*}
L_{\rm D} = \left(\frac{k_{\rm B}T\eps_0\eps_{\rm c}}{e^2N_{\rm A}}\right)^{1/2},\ \
E_0 = \left(\frac{k_{\rm B}TN_{\rm A}}{\eps_0\eps_{\rm c}}\right)^{1/2},
\textrm{ \ and \ } I_0 = \frac{k_{\rm B}TN_{\rm A}}{\eps_0\eps_{\rm c}}.
\end{equation*} 

\subsection{The Zozulya--Anderson model}
\label{Subsec_ZA}

Zozulya--Anderson model\cite{zozulya-anderson} is obtained using the above 
characteristic values and for a specific material (\chemical{LiNbO_3}) which 
imposes certain symmetries. The adiabatic assumption allows to get rid of the 
time-dependence and an asymptotic formal analysis which accounts for the 
very large donors density $N_{\rm D}$ ends the derivation. 

Dimensionless equations are obtained using $n_0$ and $N_{\rm A}$ for electron 
and ion densities respectively, $I_0$ for intensities, $E_0$ for fields and 
$\eps_{\rm c}$ for the permittivity tensor. Coefficient $\beta$ is normalized 
as a dark intensity $I_{\rm d}=\beta/sI_0$. We keep all the other notations but 
they now denote the normalized variables. The total intensity is 
$I=I_{\rm em}+I_{\rm d}$. We assume that the space charge field $\bfE_{\rm sc}$ 
derives from a potential: $L_{\rm D}\nabla\ph=-\bfE_{\rm sc}$. In the adiabatic 
assumption matter equations reduce to
\begin{eqnarray*}
& \ds I \frac{1-N_{\rm D}^+ N_{\rm A}/N_{\rm D}}{1-N_{\rm A}/N_{\rm D}} 
   = n_{\rm e} N_{\rm D}^+, \\[2mm]
& \ds L_{\rm D} \nabla \cdot \left\{n_{\rm e} \bfE_{\rm tot} 
   + L_{\rm D} \nabla n_{\rm e} 
   + \bfE_{\rm ph} I_{\rm em} \frac{1-N_{\rm D}^+ N_{\rm A}/N_{\rm D}}
                                   {1-N_{\rm A}/N_{\rm D}}\right\} = 0, \\[2mm]
& \ds -L_{\rm D}^2 \nabla \cdot (\tens\eps \nabla\ph) 
   = N_{\rm D}^+ -1 - \frac{n_0}{N_{\rm A}}n_{\rm e}.
\end{eqnarray*}
In \chemical{LiNbO_3}, $N_{\rm A}/N_{\rm D}\sim10^{-3}$ and 
$n_0/N_{\rm D}\sim10^{-6}$ and we neglect them. Finally, we make different 
assumptions on the fields: first the beam is not too thin, the photogalvanic 
and the external applied fields are not too large and therefore we may neglect 
$-L_{\rm D}^2 \nabla \cdot (\tens\eps \nabla\ph)$; second the propagation field 
amplitude is relatively small and we assimilate $\bfE_{\rm tot}$ and 
$\bfE_{\rm sc}$. This implies $n_{\rm e}=I$ and $N_{\rm D}^+=1$ and we 
have only one matter equation, namely
\begin{equation*}
\nabla I \cdot \nabla\ph + I\nabla^2\ph - \nabla^2I 
- k_{\rm D} E_{\rm ph} \bfc \cdot \nabla I = 0,
\end{equation*}
where $k_{\rm D}=1/L_{\rm D}$. To obtain a "simpler" equation, in physics 
papers the variable $U=\ph-\ln I$ is often used. This variable seems however 
to lack physical meaning. 

\noindent
The final matter equation is
\begin{equation}
\label{Eq_U} %2.7
\nabla U \cdot \nabla\ph + \nabla^2U - k_{\rm D} E_{\rm ph} \bfc \cdot \nabla I = 0.
\end{equation}

We now fix different space directions. Propagation is supposed to take place 
in the $z$-direction and $\bfk=k\bfe_z$. The two transverse directions are 
therefore $x$ and $y$. The $\bfe_x$ direction is chosen as both $\bfc$ and 
$\bfe$. If an external field is applied, it will be along $\bfe_x$ as well.
In the matter equation \eqref{Eq_U}, the quantity $\bfc \cdot \nabla I$ simply
reads $\partial_x I$. In \chemical{LiNbO_3}, $r=r_{xxx}$ is responsible
for the change of refractive index (it is $r_{xxy}$ in some other
materials) and we approximate $\tens\eps$ by $n^2$ in the expression
for $\delta n$ which becomes 
$\delta n = \frac12 n^3 r E_0 L_{\rm D}\partial_x \ph$. 
Together with Eq. \eqref{Eq_Envelope} the envelope equation now reads 
\begin{equation*}
\left[\partial_z + \frac{i}{2k} \nabla_\perp^2\right]A(\bfx) 
= - i\frac k2 n^2 r E_0 L_{\rm D}\partial_x \ph A(\bfx).
\end{equation*}

The last step is to have dimensionless space variables. We set 
$\alpha = \frac k2 n^2 r E_0$ which has the dimension of the inverse of a 
space variable. We denote $z'=|\alpha|z$, $(x',y')=\sqrt{k|\alpha|}(x,y)$, 
$A'= A/\sqrt{I_0I_{\rm d}}$, $\ph'=\sqrt{k|\alpha|}\ph/k_{\rm D}$ and 
$U'=\sqrt{k|\alpha|}U/k_{\rm D}$. 
The last approximations are now $U'=\ph'$, $k\gg|\alpha|$ and 
$E_{\rm ph}\sim E_0$, and omitting primes:
\begin{eqnarray*}
\left[\partial_z-\frac i2 \nabla_\perp^2\right]A & = & -iA\partial_x\ph, \\[2mm]
\nabla_\perp^2 \ph + \nabla_\perp \ln(1+|A|^2)\cdot \nabla_\perp \ph 
& = & \partial_x\ln(1+|A|^2).
\end{eqnarray*}
These equations are usually referred to as a model derived in
\refcite{zozulya-anderson} but only seeds of these equations are derived 
there usually including many other terms and especially time derivatives.

In the wide literature devoted to photorefractive media, many equations are 
written which resemble those above but with different choices of asymptotic 
approx\-imations. In particular numerical results are very often obtained keeping 
the time in the matter equations (see \refcite{stepken-kaiser-belic-krolikowski} or
\refcite{wolfersberger-fressengeas-maufroy-kugel}).

\subsection{Mathematical setting}

If we look at a wider class of materials we may have different signs for the 
nonlinearity (in reference to the cubic nonlinear Schr\"odinger equation, 
the case $a=1$ is classically called the focusing case, and $a=-1$ the 
defocusing case). Besides mathematicians are more accustomed to use $t$ as the 
evolution variable. We will therefore consider the system
\begin{equation}
\label{Eq_Aln} %2.8
\left\{
\begin{array}{ll}
i\partial_t A + \Delta A = -aA\partial_x\ph,              & a=\pm1, \\
\Delta \ph + \nabla \ln(1+|A|^2) \cdot \nabla \ph  = \partial_x\ln(1+|A|^2),
\end{array}
\right.
\end{equation}
where $\Delta = \partial_x^2+\partial_y^2$ or $\Delta=\partial_x^2$. 

These expressions with logarithms are widely used in the physics literature, 
maybe because they are the starting point of solitonic studies and logarithms 
appear
naturally in the expression of solitary waves (see Sec. \ref{Sec_1D}). This 
form is however cumbersome to handle for the mathematical analysis, and it is 
much more convenient to cast \eqref{Eq_Aln} as
\begin{equation}
\label{Eq_Adiv} %2.9
\left\{
\begin{array}{l}
i\partial_t A + \Delta A = -aA\partial_x\ph, \\[2mm]
\div \left((1+|A|^2) \nabla \ph\right) = \partial_x(|A|^2),
\end{array}
\right.
\end{equation}
which is closer to the original Kukhtarev equations.

We have seen that the main effects take place in the $t$- (propagation) and the 
$x$-directions (drift, anisotropic diffusion, external field, polarization). 
It is therefore natural to study the equations with no dependence in the $y$ 
variable. In the one-dimensional case, we infer immediately from the last 
equation in System \eqref{Eq_Adiv} that 
$(1+|A|^2) \partial_x \ph = |A|^2 - C(t)$ where the constant $C(t)$ is given by 
the boundary conditions. If no external field is applied $C(t)\equiv0$. This is 
the case for bright solitary waves (see \refcite{mamaev-saffman-zozulya2}). 
In the case of dark solitary waves $C(t) = \lim_{x\to\pm\infty} |A|^2$ 
(see \refcite{mamaev-saffman-zozulya1}), which does not depend on $t$ 
either. In both cases, System \eqref{Eq_Adiv} reduces to the saturated NLS 
equation
\begin{equation}
\label{Eq_SaturatedNLS} %2.10
i \partial_t A + \partial_x^2 A = -a \frac{|A|^2 - |A_\infty|^2}{1+|A|^2}A.
\end{equation}
In the sequel we will mainly consider the case when $A_\infty=0$ and show 
that, in some sense, the dynamics of \eqref{Eq_Adiv} is similar to that of 
\eqref{Eq_SaturatedNLS} which we will recall in Sec. \ref{Sec_saturated}.

In the two-dimensional case \eqref{Eq_Adiv} can be viewed as a saturated version 
of a Davey--Stewartson system. Namely, replacing $1+|A|^2$ by 1 in the L.H.S. 
of \eqref{Eq_Adiv} yields
\begin{equation*}
\left\{
\begin{array}{rcl}
i\partial_t A + \Delta A & = & -a A\partial_x\ph, \\[2mm]
\Delta \ph & = & \partial_x(|A|^2),
\end{array}
\right.
\end{equation*}
which is the Davey--Stewartson system of the elliptic--elliptic type (see 
\refcite{ghidaglia-saut1}). 

\section{The Saturated NLS Equation}
\label{Sec_saturated}

We review here some mathematical facts, more or less known, on the saturated 
NLS equation
\begin{equation}
\label{Eq_NLSCauchy} %3.1
\left\{
\begin{array}{ll}
\ds i\partial_t A + \Delta A = - a\frac{|A|^2A}{1+|A|^2},\ & a=\pm1, \\[2mm]
A(\bfx,0) = A_0(\bfx),
\end{array}
\right.
\end{equation}
where $A=A(\bfx,t)$ and $\bfx\in\bbR^d$. We have derived this equation for 
$d=1$, but give here results for a general $d$. This equation is also derived
in other contexts, for example the propagation of a laser beam in gas vapors
\cite{tikhonenko-christou-luther-davies}.

\subsection{The Cauchy problem}

The Cauchy problem \eqref{Eq_NLSCauchy} can be solved in $L^2$ and in the energy 
space $H^1$.

\begin{theorem}
\label{Th_NLSCauchy}
{\em (i)} Let $A_0\in L^2(\bbR^d)$. Then there exists a unique solution 
$A\in\calC(\bbR;L^2(\bbR^d))$ of \eqref{Eq_NLSCauchy} which satisfies 
furthermore
\begin{equation}
\label{Eq_NLSL2Claw}
\int_{\bbR^d} |A(\bfx,t)|^2 d\bfx = \int_{\bbR^d} |A_0(\bfx)|^2 d\bfx, 
\hspace{1cm} t\in\bbR. 
\end{equation}
{\em (ii)} Let $A_0\in H^1(\bbR^d)$. Then the solution above satisfies 
$A\in\calC(\bbR;H^1(\bbR^d))$ and
\begin{eqnarray}
&& \int_{\bbR^d} 
\left[|\nabla A(\bfx,t)|^2 d\bfx + a \ln (1+|A(\bfx,t)|^2) \right]
d\bfx \nonumber \\
\label{Eq_NLSH1Claw}
&& \hspace{1cm} = \int_{\bbR^d} 
\left[|\nabla A_0(\bfx)|^2 d\bfx + a \ln (1+|A_0(\bfx)|^2) \right]
d\bfx, \hspace{1cm} t\in\bbR. 
\end{eqnarray}
\end{theorem}

\proof
The norm conservations \eqref{Eq_NLSL2Claw} and \eqref{Eq_NLSH1Claw} result 
from multiplying \eqref{Eq_NLSCauchy} by $\bar A$ and $\partial_t \bar A$ 
respectively and integrating the complex and real parts respectively. This 
formal proof is justified by the standard truncation process. 

Let $S(t)$ be the group operator associated to the linear Schr\"odinger equation 
$i\partial_t A + \Delta A =0$. Then the Duhamel formula for \eqref{Eq_NLSCauchy} 
reads
\begin{equation}
\label{Eq_3.4}
A(\bfx,t) = S(t) A_0(\bfx) 
- a \int_0^t S(t-s) \frac{|A(\bfx,s)|^2}{1+|A(\bfx,s)|^2}A(\bfx,s)\ ds.
\end{equation}
Since $x\mapsto x/(1+x)$ is Lipschitz, we easily infer that the R.H.S. of 
\eqref{Eq_3.4} defines a contraction on a suitable ball of 
$\calC([0,T];L^2(\bbR^d))$ for some $T>0$. The local well-posedness in 
$L^2(\bbR^d))$ follows. Global well-posedness is derived from the conservation
law \eqref{Eq_NLSL2Claw}.

The $H^1$ theory follows the same argument, noticing that
\begin{equation*}
\left|\nabla\left(\frac{|A|^2}{1+|A|^2}A\right)\right| 
= \left|\frac{A^2}{(1+|A|^2)^2}\nabla \bar A 
        + \frac{|A|^2}{(1+|A|^2)^2}\nabla A\right|
\leq \frac12|\nabla A|.
\end{equation*}
\endproof

\remark As a consequence of \eqref{Eq_NLSL2Claw}, \eqref{Eq_NLSH1Claw} and 
$\ln(1+|A|^2)\leq|A|^2$, we obtain the uniform bound 
\begin{equation}
\label{Eq_NLSH1Est}
\int_{\bbR^d} |\nabla A(\bfx,t)|^2 d\bfx 
\leq \int_{\bbR^d} |A_0(\bfx)|^2 d\bfx 
+ \int_{\bbR^d} |\nabla A_0(\bfx)|^2 d\bfx, 
\hspace{1cm} t\in\bbR. 
\end{equation}
Contrarily to the context of the usual nonlinear cubic Schr\"odinger equation,
this bound does not depend on the sign of $a$ and in particular saturation 
means that no blow-up is occurs. 

\subsection{Solitary waves --- one-dimensional results}
\label{Sec_1D}

In the one-dimensional case, it is possible to compute first integral 
formulations of the solitary waves. \\

Bright solitary waves are sought for in the form $A(x,t) = e^{i\omega t}u(x)$ 
(see \refcite{mamaev-saffman-zozulya2}), where $A$ is a solution to 
\eqref{Eq_SaturatedNLS} with $A_\infty=0$. The function $u$ is supposed to have
a maximum at $x=0$ ($u(0)=u_m>0$ and $u'(0)=0$), therefore 
\begin{equation*}
[u'(x)]^2 = (\omega-a)[u^2(x)-u_m^2] + a [\ln(1+u^2)-\ln(1+u_m^2)].
\end{equation*}
We furthermore want that for $x\to\infty$, $u(x)\to0$ and $u'(x)\to0$. This
yields a unique possible frequency for the solitary wave, namely
\begin{equation*}
\omega = a\left(1-\frac{\ln(1+u_m^2)}{u_m^2}\right)
\end{equation*}
and
\begin{equation*}
[u'(x)]^2 = a \left(-\frac{u^2(x)}{u_m^2}\ln(1+u_m^2)+\ln(1+u^2)\right).
\end{equation*}
Since $u_m$ is supposed to be the maximum of $u$, this quantity is positive 
only if $a=1$ (focusing case) and the bright soliton is solution to the 
first order equation:
\begin{equation*}
u'(x) = -\textrm{sign}(x) \sqrt{\ln(1+u^2)-\frac{u^2}{u_m^2}\ln(1+u_m^2)} 
\textrm{ with } \omega = 1-\frac{\ln(1+u_m^2)}{u_m^2}.
\end{equation*}

Dark solitary waves are sought for in the form $A(x,t) = u(x)$ (see 
\refcite{mamaev-saffman-zozulya1}) where $A$ is solution to
\eqref{Eq_SaturatedNLS}. There is no time-dependence. We assume that 
$\lim_{x\to\pm\infty}u'(x)=0$ and consistently with $A_\infty\neq0$,  
\begin{equation*}
\lim_{x\to+\infty}u(x)=-\lim_{x\to-\infty}u(x)=u_\infty.
\end{equation*}
Then 
\begin{equation*}
[u'(x)]^2 
= a \left( - (u^2-u_\infty^2) 
+ (1+u_\infty^2)\ln\left(\frac{1+u^2}{1+u_\infty^2}\right) \right).
\end{equation*}
At the origin $u(0)=0$ and we want more generally that $|u(x)|\leq|u_\infty|$.
Therefore, dark solitary waves only exist if $a=-1$ (defocusing case). In this 
context $u(x)$ is a monotonous function and is solution to the first order 
equation:
\begin{equation*}
u'(x) = \textrm{sign} (u_\infty) \sqrt{u^2-u_\infty^2 
- (1+u_\infty^2)\ln\left(\frac{1+u^2}{1+u_\infty^2}\right)}.
\end{equation*}
For both bright and dark solitary waves, no explicit solution is known.

\subsection{Solitary waves --- \textit{a priori} estimates and non existence}

Consider now the solitary wave solutions of \eqref{Eq_NLSCauchy} in any dimension 
$d$, that is solutions of the type $A(\bfx,t)=e^{i\omega t}U(\bfx)$, where 
$U\in H^1(\bbR^d)$ (we thus are only concerned with "bright" solitary waves). 
A solitary wave is a solution of the elliptic equation
\begin{equation}
\label{Eq_NLSSW}
-\Delta U + \omega U = a \frac{|U|^2U}{1+|U|^2}, \hspace{1cm} U\in H^1(\bbR^d).
\end{equation}
A trivial solution is $U\equiv0$. We seek other nontrivial solutions.

\begin{lemma}
Any $H^1(\bbR^d)$ solitary wave satisfies
\begin{equation}
\label{Eq_NLSSWEnergy}
\int_{\bbR^d} \left[ |\nabla U|^2 
+ \left( \omega - a \frac{|U|^2}{1+|U|^2} \right) |U|^2 \right] d\bfx =0
\end{equation}
(energy identity)
\begin{equation}
\label{Eq_NLSSWPohozaev}
(d-2) \int_{\bbR^d} |\nabla U|^2 d\bfx + d\omega \int_{\bbR^d} |U|^2 d\bfx 
- ad  \int_{\bbR^d} \left[ |U|^2 - \ln (1+|U|^2) \right] d\bfx =0
\end{equation}
(Pohozaev identity).
\end{lemma}

\proof
As for Theorem \ref{Th_NLSCauchy}, \eqref{Eq_NLSSWEnergy} results from 
multiplying \eqref{Eq_NLSSW} by $\bar U$ and inte\-grating. To get 
\eqref{Eq_NLSSWPohozaev}, one multiplies \eqref{Eq_NLSSW} by 
$x_k\partial_{x_k}\bar U_k$, integrates the real part, and sums from 1 to $d$. 
This is justified by a standard truncation argument.
\endproof

\begin{corollary}
\label{Cor_NLSnoSW}
No nontrivial solitary wave (solution of \eqref{Eq_NLSSW}) exists when
\setlength\leftmargini{2pc}
\begin{itemize}
\item[{\em (i)}] $a=-1$ {\em(}defocusing case{\em)}, for $\omega\geq0$.
\item[{\em (ii)}] $a=1$ {\em(}focusing case{\em)} and $\omega\geq1$.
\item[{\em (iii)}] $a=\pm1$ if $\omega<0$ provided 
$|U|^2/(1+|U|^2) = O(1/|\bfx|^{1+\eps})$, $\eps>0$ as $|\bfx|\to+\infty$.
\end{itemize}
\end{corollary}

\proof
Identity \eqref{Eq_NLSSWEnergy} implies that no solitary wave may exist when 
$a=-1$ and $\omega\geq0$ or $a=1$ and $\omega\geq1$. When $d=1,2$, Eq. 
\eqref{Eq_NLSSWPohozaev} implies that no solitary wave exist when $\omega\leq0$ 
and $a=1$. Recall $d=2$ is the physical case. The remaining cases ($\omega<0$, 
$a=-1$ or $a=1$, $d\geq3$) follow from the classical result of Kato\cite{kato} 
on the absence of embedded eigenvalues. Indeed, we can write \eqref{Eq_NLSSW} as
\begin{equation*}
\Delta U + (-\omega + V(\bfx)) U = 0, \hspace{1cm} V(\bfx)=a\frac{|U|^2}{1+|U|^2},
\end{equation*}
assuming furthermore that $V(\bfx)=O(1/|\bfx|^{1+\eps})$, $\eps>0$, as 
$|\bfx|\to\infty$. A proof for $d=3,4$ or $d\geq5$, $\omega\leq-(d-2)/2$ with no
decaying assumption is given in Appendix.
\endproof

\begin{corollary}
\label{Cor_NLSSWExist}
Solitary waves may exist only when $a=1$ and $0<\omega<1$. 
\end{corollary}

Corollary \ref{Cor_NLSSWExist} is consistent with the one-dimensional "explicit" 
result. We first have a classical regularity and decay result. 

\begin{proposition}
\label{Prop_NLSSWDecay}
Let $a=1$ and $0<\omega<1$. Then any $U\in H^1(\bbR^d)$ solution of 
\eqref{Eq_NLSSW} satisfies
\begin{equation}
\label{Eq_NLSSWHinfty}
U \in H^\infty(\bbR^d),
\end{equation}
\begin{equation}
\label{Eq_NLSSWDecay}
e^{\delta|\bfx|} U \in L^\infty(\bbR^d) \textrm{ for any } \delta<\omega/2.
\end{equation}
\end{proposition}

\proof
$U\in H^\infty(\bbR^d)$ results trivially from a bootstrapping argument using 
$|U|^2/(1+|U|^2)<1$. To prove \eqref{Eq_NLSSWDecay}, we first derive the estimate
\begin{equation}
\label{Eq_NLSSWH1FiniteWeight}
\int_{\bbR^d} e^{\omega|\bfx|} \left[ |\nabla U|^2 + |U|^2 \right] d\bfx 
< +\infty.
\end{equation}
In fact, as in Cazenave\cite{cazenave1,cazenave2}, we multiply \eqref{Eq_NLSSW} 
by $e^{\omega|\bfx|}\bar U$ and integrate the real part (this formal argument is 
made rigorous by replacing $e^{\omega|\bfx|}$ by 
$e^{\omega|\bfx|/(1+\eps|\bfx|)}$, $\eps>0$, $\eps\to0$) to get
\begin{equation}
\label{Eq_NLSSWH1Weight}
\int_{\bbR^d} e^{\omega|\bfx|} \left[ |\nabla U|^2 + \omega |U|^2 \right] d\bfx
\leq \int_{\bbR^d} e^{\omega|\bfx|} \frac{|U|^4}{1+|U|^2} d\bfx 
+ \omega \int_{\bbR^d} e^{\omega|\bfx|} |U||\nabla U| d\bfx.  
\end{equation}
By \eqref{Eq_NLSSWHinfty} there exists $R>0$ such that $|U|^2/(1+|U|^2)<\omega/4$ on 
$\bbR^d\setminus B_R$. Thus we infer from \eqref{Eq_NLSSWH1Weight} that
\begin{eqnarray*}
\int_{\bbR^d} e^{\omega|\bfx|} \left[ |\nabla U|^2 + \omega |U|^2 \right] d\bfx
& \leq & \int_{B_R} e^{\omega|\bfx|} \frac{|U|^4}{1+|U|^2} d\bfx 
+ \frac\omega4 \int_{\bbR^d} e^{\omega|\bfx|} |U|^2 d\bfx \\ 
&& + \frac\omega2 \left(\int_{\bbR^d}  e^{\omega|\bfx|}|U|^2 d\bfx
+ \int_{\bbR^d}  e^{\omega|\bfx|}|\nabla U|^2 d\bfx \right) 
\end{eqnarray*}
which implies \eqref{Eq_NLSSWH1FiniteWeight}.

Now we write $U$ as a convolution
\begin{equation}
\label{Eq_NLSSWCovolution}
U(\bfx) = H_\omega \star \frac{|U|^2U}{1+|U|^2},
\hspace{1cm}\textrm{ where }
H_\omega = \calF^{-1} \left(\frac1{\omega+|\xi|^2}\right).
\end{equation}
As it is well known (\refcite{abramowitz-steglun}), 
$H_\omega(\bfx)=\omega^{(d-2)/2} G_1(\omega^{1/2}\bfx)$ where
\begin{equation*}
G_1(z) = 
\begin{cases}
C|z|^{(2-d)/2} K_{(d-2)/2} (|z|), & d\geq3, \\[2mm]
K_0 (|z|), & d=2,
\end{cases}
\end{equation*}
where $K_\nu$ is the modified Bessel function of order $\nu$. Furthermore (see 
\refcite{abramowitz-steglun}), one has the asymptotic behavior:
\begin{equation}
\label{Eq_Bessel}
\left\{
\begin{array}{rcll}
K_\nu (z) & \sim & \frac12 \Gamma(\nu) \left(\frac12|z|\right)^{-\nu}, 
& \textrm{ for } \nu>0, \textrm{ as } |z|\to0, \\[2mm]
K_0 (z) & \sim & - \ln(|z|), & \textrm{ as } |z|\to0, \\[2mm]
K_\nu (z) & \sim & C|z|^{-1/2} e^{-|z|}, & \textrm{ as } |z|\to\infty.
\end{array}
\right.
\end{equation}
We infer from \eqref{Eq_NLSSWCovolution} that
\begin{equation}
\label{Eq_NLSSWCovolutionEst}
e^{\delta|\bfx|}|U(\bfx)| 
\leq \int_{\bbR^d} e^{\delta|\bfx-\bfx'|}H_\omega(\bfx-\bfx')e^{\delta|\bfx'|} 
\frac{|U|^3}{1+|U|^2}(\bfx')d\bfx'.
\end{equation}
Since by \eqref{Eq_Bessel} $e^{\delta|\bfx|}H_\omega(\bfx)\in L^2(\bbR^d)$ for
$0<\delta<\omega^{1/2}$, and by \eqref{Eq_NLSSWH1FiniteWeight} 
$e^{\delta|\bfx|} |U|^3/(1+|U|^2) \in L^2(\bbR^d)$ for 
$\delta\leq\omega/2<\omega^{1/2}$, we deduce from \eqref{Eq_NLSSWCovolutionEst} 
that $e^{\delta|x|}U\in L^\infty(\bbR^d)$ for $0<\delta<\omega/2$.
\endproof

\remark 
Actually, the saturated cubic NLS equation should involve a small 
parameter $\eps>0$, namely, in the focusing case, we should consider instead 
of \eqref{Eq_NLSCauchy} 
\begin{equation}
\label{Eq_NLSeps}
i\partial_t A^\eps + \Delta A^\eps = - \frac{|A^\eps|^2A^\eps}{1+\eps|A^\eps|^2}.
\end{equation}
Theorem \ref{Th_NLSCauchy} is of course still valid for a fixed $\eps>0$, but 
\eqref{Eq_NLSH1Claw} and \eqref{Eq_NLSH1Est} should be replaced by
\begin{equation*}
\begin{aligned}
\int_{\bbR^d} &
\left[|\nabla A(\bfx,t)|^2 d\bfx + \frac1{\eps^2} \ln (1+|A(\bfx,t)|^2 \right]
d\bfx \\
& = \int_{\bbR^d} 
\left[|\nabla A_0(\bfx)|^2 d\bfx + \frac1{\eps^2} \ln (1+|A_0(\bfx)|^2 \right]
d\bfx,
\end{aligned}
\end{equation*}
\begin{equation*}
\int_{\bbR^d} |\nabla A(\bfx,t)|^2 d\bfx 
\leq \frac1\eps \int_{\bbR^d} |A_0(\bfx)|^2 d\bfx 
+ \int_{\bbR^d} |\nabla A_0(\bfx)|^2 d\bfx. 
\end{equation*}
For solitary waves $A^\eps(\bfx,t)=e^{i\omega t} U(\bfx)$, \eqref{Eq_NLSeps} 
reduces to the elliptic equation
\begin{equation*}
-\Delta U + \omega U = \frac{|U|^2U}{1+\eps|U|^2}.
\end{equation*}
Setting $V=\eps^{1/2}U$, one obtains
\begin{equation*}
-\Delta V + \omega V = \frac1\eps\ \frac{|V|^2V}{1+|V|^2}.
\end{equation*}
The only possible range for the existence of nontrivial $H^1$ solitary waves
is $\omega\in]0,1/\eps[$. Proposition \ref{Prop_NLSSWDecay} is still valid for 
$\omega$ in this range.

\subsection{Solitary waves --- existence results}

We now turn to the existence of non-trivial $H^2$ solutions of
\begin{equation*}
-\Delta U + \omega U = \frac{|U|^2U}{1+|U|^2}
\end{equation*}
when $0<\omega<1$. We will look for real radial solutions
$U(\bfx)=u(|\bfx|)\equiv u(r)$ and thus consider 
the ODE problem
\begin{equation}
\label{Eq_NLSSWRadial}
\left\{
\begin{array}{l}
\ds - u'' - \frac{d-1}{r} u' + \omega u = \frac{u^3}{1+u^2}, \\
u \in H^2(]0,\infty[), \hspace{1cm} u'(0)=0.
\end{array}
\right.
\end{equation}
We recall a classical result of Berestycki \textit{et al.}% 
\cite{berestycki-lions-peletier}

\begin{theorem}[\refcite{berestycki-lions-peletier}, p. 143]
\label{Th_BLP}
Let $g$ be a locally Lipschitz continuous function on $\bbR_+=[0,+\infty[$ 
such that $g(0)=0$, satisfying the following hypotheses.
\setlength\leftmargini{2pc}
\begin{itemize}
\item[{\em (H1)}] $\alpha=\inf\{\zeta>0,\ g(\zeta) \geq 0\}$ exists and $\alpha>0$.
\item[{\em (H2)}] Let $G(t) = \int_0^t g(s) ds$. There exists $\zeta>0$ such that
$G(\zeta)>0$. \\
Let $\zeta_0=\inf\{\zeta>0,\ G(\zeta) \geq 0\}$. In view of {\em (H1)} and 
{\em (H2),} $\zeta_0$ exists and $\zeta_0>\alpha$. 
\item[{\em (H3)}] $\lim_{s\searrow\alpha} g(s)/(s-\alpha)>0$
\item[{\em (H4)}] $g(s)>0$ for $s\in]\alpha,\zeta_0]$. \\ 
Let $\beta=\inf\{\zeta>\zeta_0,\ G(\zeta) \geq 0\}$. In view of {\em (H4),} 
$\zeta_0<\beta\leq+\infty$.
\item[{\em (H5)}] If $\beta=+\infty$, then $\lim_{s\to+\infty} g(s)/s^l=0$ 
for some $l<(d+2)/(d-2)$ (if $d=2$, we may choose for $l$ any finite real 
number).
\end{itemize}
Let us consider the Cauchy problem
\begin{equation}
\label{Eq_BPL}
\left\{
\begin{array}{ll}
\ds - u'' - \frac{d-1}{r} u' = g(u),\ & r>0, \\
u(0)=\zeta, & u'(0)=0.
\end{array}
\right.
\end{equation}
Then there exists $\zeta\in]\zeta_0,\beta[$ such that \eqref{Eq_BPL} has a 
unique solution satisfying $u(r)>0$ for $r\in]0,+\infty[$, $u'(r)<0$ for 
$r\in]0,+\infty[$ and $\lim_{r\to+\infty}u(r)=0$. If in addition 
$\limsup_{s\searrow0} g(s)/s<0$, then there exists $C>0$ and $\delta>0$ 
such that $0<u(r)\leq Ce^{-\delta r}$, for $0\leq r<+\infty$.
\end{theorem}

\begin{theorem}
If $a=1$ and $0<\omega<1$, there exists a nontrivial positive solution of 
\eqref{Eq_NLSSWRadial}.
\end{theorem}

\begin{figure}[b]
\centerline{\psfig{file=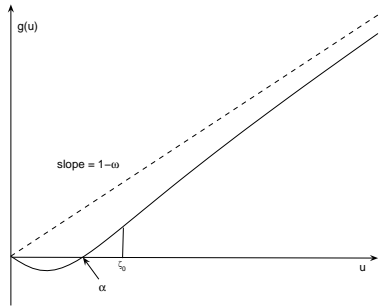,width=5cm}}
\vspace*{8pt}
\caption{\label{Fig_1}Graph of function $g$}
\end{figure}

\proof
The case $d=1$ has been addressed in Sec. \ref{Sec_1D}. Consider now 
$d\geq2$. We apply Theorem \ref{Th_BLP} with 
\begin{equation*}
g(u)=-\omega u + \frac{u^3}{1+u^2},
\end{equation*} 
which graph is displayed in Fig. \ref{Fig_1}.

Note that $\alpha=\sqrt{\omega/(1-\omega)}$ which yields (H1). Setting 
$u=\sqrt{\omega/(1-\omega)}+\eps$, one easily checks that
\begin{equation*}
\frac{g(u)}{u-\alpha} = 2\omega(1-\omega) + O(\eps),
\end{equation*}
and (H3) is satisfied. One computes $G(u)=(1-\omega)u^2-\frac12 \ln(1+u^2)$, 
which obviously satisfies (H2) and (H4) with $\beta=+\infty$. Last (H5) holds 
true (for $l>1$).
\endproof

\remark 
$u$ satisfies the decay rate of Proposition \ref{Prop_NLSSWDecay}.

\section{The Zozulya--Anderson System}
\label{Sec_Cauchy}

\subsection{Estimate on the potential}

We now restrict to the space-dimension $d=2$ which is the context of the
derivation. To mimic the proof for the Cauchy problem in the one-dimensional 
case, we would like to express $\ph$ in terms of $A$ for say 
$A\in L^2(\bbR^2)$. With such a data $A$, we indeed have a unique $\ph$ in some 
convenient space but no Lipschitz regularity for the mapping $A\mapsto\ph$,
which is required to perform some fixed point procedure. To ensure this we 
will have to assume $A\in H^2(\bbR^2)$. 

To derive the first estimates, we consider time as a parameter and do not 
express it. We therefore introduce the weighted homogeneous Sobolev space 
\begin{equation*}
H = \{\ph\in\calS'(\bbR^d),\ (1+|A|^2)^{1/2}\nabla\ph\in L^2(\bbR^d)\}\slash\bbR
\end{equation*}
together with its natural Hilbertian structure.

\begin{lemma}
\label{Lm_Invdivgrad}
{\em (i)} Let $A\in L^2(\bbR^2)$. There exists a unique $\ph\in H$ solution of
\begin{equation}
\label{Eq_ZAPotential}
\div ((1+|A|^2)\nabla\ph) = \partial_x (|A|^2) \textrm{ in } \calD'(\bbR^2)
\end{equation}
such that
\begin{equation}
\label{Eq_ZAH1Est}
\int_{\bbR^2} (1+\frac12|A|^2)|\nabla\ph|^2 d\bfx 
\leq \frac12 \int_{\bbR^2} |A|^2 d\bfx. 
\end{equation}
{\em (ii)} If furthermore $A\in H^2(\bbR^2)$, then $\nabla\ph\in H^2(\bbR^2)$ and 
there exists a polynomial $P$ vanishing at 0 such that
\begin{equation}
\label{Eq_P}
\|\nabla \ph\|_{H^2(\bbR^2)} \leq P(\|A\|_{H^2(\bbR^2)}).
\end{equation}
\end{lemma}

\proof
(i) We define a smoothing sequence $(\theta_\eps)_{\eps>0}$ with
$\int_{\bbR^2}\theta_\eps d\bfx=1$ and $A_\eps=A\star\theta_\eps$ is such that 
$A_\eps\to A\in L^2(\bbR^2)$. In particular
\begin{equation}
\label{Eq_ZAUnif}
\|A_\eps\|_{L^2(\bbR^2)} \leq \|A\|_{L^2(\bbR^2)}. 
\end{equation}
By Riesz theorem there exists a unique solution to 
\begin{equation}
\label{Eq_ZAPotentialeps}
\div ((1+|A_\eps|^2)\nabla\ph_\eps) = \partial_x (|A_\eps|^2),
\end{equation}
i.e. 
\begin{equation}
\label{Eq_ZAEllipticeps}
-\Delta \ph_\eps - \div(|A_\eps|^2\nabla\ph_\eps) = \partial_x (|A_\eps|^2)
\end{equation}
after noticing that the R.H.S. of Eq. \eqref{Eq_ZAPotentialeps} defines a 
linear continuous form on $H$ given by
\begin{equation*}
\langle \partial_x(|A_\eps|^2),\psi \rangle 
= \int_{\bbR^2} |A_\eps|^2 \partial_x\psi d\bfx.
\end{equation*}

\noindent
Now we get from \eqref{Eq_ZAPotentialeps}
\begin{eqnarray*}
\int_{\bbR^2} (1+|A_\eps|^2)|\nabla\ph_\eps|^2 d\bfx 
& = & - \int_{\bbR^2} |A_\eps|^2 \partial_x\ph_\eps d\bfx \\
& \leq & \frac12 \int_{\bbR^2} |A_\eps|^2  d\bfx 
+ \frac12  \int_{\bbR^2} |A_\eps|^2 (\partial_x\ph_\eps)^2 d\bfx,
\end{eqnarray*}
which yields (together with \eqref{Eq_ZAUnif})
\begin{equation}
\label{Eq_ZAH1Esteps}
\int_{\bbR^2} (1+\frac12|A_\eps|^2)|\nabla\ph_\eps|^2 d\bfx 
\leq \frac12 \int_{\bbR^2} |A_\eps|^2 d\bfx
\leq \frac12 \int_{\bbR^2} |A|^2 d\bfx. 
\end{equation}
Up to the extraction of a sub-sequence, we have $\nabla \ph_\eps\to\nabla\ph$
and $\partial_x (|A_\eps|^2)\to\partial_x (|A|^2)$ in $\calD'(\bbR^2)$. From 
Eq. \eqref{Eq_ZAH1Esteps}, $A_\eps\nabla\ph_\eps\rightharpoonup B$ weakly
in $L^2$ and for all $\psi\in\calD$, 
\begin{equation*}
\int_{\bbR^2} A_\eps\nabla\ph_\eps\cdot\nabla\psi d\bfx 
\to
\int_{\bbR^2} A\nabla\ph\cdot\nabla\psi d\bfx, 
\end{equation*}
therefore $B=A\nabla\ph$.  Since
$||A_\eps|^2\nabla\ph_\eps|=|A_\eps||A_\eps\nabla\ph_\eps|$,
$|A_\eps|^2\nabla\ph_\eps\to|A|^2\nabla\ph$ in $\calD'$. We can pass to the
limit in Eq. \eqref{Eq_ZAEllipticeps} and obtain 
\begin{equation*}
-\Delta \ph - \div(|A|^2\nabla\ph) = \partial_x (|A|^2) 
\textrm{ in } \calD'(\bbR^2),
\end{equation*}
i.e. $\div ((1+|A|^2)\nabla\ph) = \partial_x (|A|^2)$ and deduce estimate 
\eqref{Eq_ZAH1Est} from \eqref{Eq_ZAH1Esteps}. This yields the existence of
$\ph\in H$. The uniqueness is straightforward: two solutions $\ph_1$ and
$\ph_2$ would satisfy
\begin{equation*}
\int_{\bbR^2} (1+\frac12|A|^2)|\nabla(\ph_1-\ph_2)|^2 d\bfx = 0, 
\textrm{ i.e. }\nabla(\ph_1-\ph_2)=0\textrm{ a.e. }
\end{equation*}
and hence be equal in $H$. \\[2mm]
(ii) We first notice that $|A|^2\Delta\ph$ is meaningful in 
$H^{-1}(\bbR^2)$. Actually, for any $\psi\in H^1(\bbR^2)$, one defines
\begin{equation*}
\langle |A|^2\Delta\ph,\psi\rangle_{H^{-1}(\bbR^2),H^1(\bbR^2)}
= \langle \Delta\ph,|A|^2\psi\rangle_{H^{-1}(\bbR^2),H^1(\bbR^2)},
\end{equation*}
which makes sense since $|A|^2\psi\in H^1(\bbR^2)$ for $A\in H^2(\bbR^2)$,
$\psi\in H^1(\bbR^2)$. Thus we can write 
\eqref{Eq_ZAPotential} as
\begin{equation*}
(1+|A|^2)\Delta\ph = - \nabla|A|^2\cdot\nabla\ph + \partial_x(|A|^2),
\end{equation*}
and
\begin{equation*}
\Delta\ph = - \frac{\nabla|A|^2}{1+|A|^2}\cdot\nabla\ph 
+ \frac{\partial_x(|A|^2)}{1+|A|^2} =: F.
\end{equation*}
We claim that $F\in L^r(\bbR^2)$, for any $r\in(1,2)$, with 
\begin{equation*}
\|F\|_{L^r(\bbR^2)} \leq C \|A\|_{L^2(\bbR^2)} \|A\|_{H^2(\bbR^2)}.
\end{equation*}
First, $|\nabla|A|^2\cdot\nabla\ph| \leq 2|\nabla A||A\nabla\ph|$ and by
H\"older
\begin{equation*}
\| \nabla|A|^2\cdot\nabla\ph \|_{L^r(\bbR^2)} 
\leq 2 \|\nabla A\|_{L^p(\bbR^2)} \|A\nabla \ph\|_{L^2(\bbR^2)}
\end{equation*}
for any $1<r<2$ and $p=2r/(2-r)\in(2,\infty)$. Since
$\|A\nabla\ph\|_{L^2(\bbR^2)} \leq \|A\|_{L^2(\bbR^2)}$ and 
$H^1(\bbR^2)\subset L^q(\bbR^2)$ for all $q>2$, we obtain that
\begin{equation*}
\left\|\frac{\nabla|A|^2}{1+|A|^2}\cdot\nabla\ph\right\|_{L^r(\bbR^2)}
\leq C \|A\|_{H^2(\bbR^2)} \|A\|_{L^2(\bbR^2)}, \hspace{1cm} 1<r<2.  
\end{equation*}
Similarly
\begin{eqnarray*}
\left\|\frac{\partial_x|A|^2}{1+|A|^2}\right\|_{L^r(\bbR^2)}
& \leq & 2 \|A \partial_x A\|_{L^r(\bbR^2)} 
\leq 2 \|A\|_{L^2(\bbR^2)} \|\partial_x A\|_{L^p(\bbR^2)} \\[2mm]
& \leq & C \|A\|_{H^2(\bbR^2)} \|A\|_{L^2(\bbR^2)}, \hspace{1cm} 1<r<2.  
\end{eqnarray*}

By elliptic regularity, we infer thus that for any $r$, $1<r<2$,
\begin{equation*}
\|\nabla\ph\|_{W^{1,r}(\bbR^2)} \leq C \|A\|_{H^2(\bbR^2)} \|A\|_{L^2(\bbR^2)}.  
\end{equation*}
By Sobolev embedding,
\begin{equation*}
\|\nabla\ph\|_{L^q(\bbR^2)} 
\leq C \|\nabla\ph\|_{W^{1,r}(\bbR^2)} 
\leq C \|A\|_{H^2(\bbR^2)} \|A\|_{L^2(\bbR^2)}.  
\end{equation*}
for $\frac 1q = \frac1r - \frac12$, i.e. $q=2r/(2-r)$ for all $r$, $1<r<2$. Thus
for any $p>2$
\begin{eqnarray*}
\left\|\frac{\nabla|A|^2}{1+|A|^2}\cdot\nabla\ph\right\|_{L^p(\bbR^2)}
& \leq & \|\nabla \ph\|_{L^{2p}(\bbR^2)} \|\nabla |A|^2\|_{L^{2p}(\bbR^2)} \\[2mm]
& \leq & C \|A\|_{H^2(\bbR^2)} \|A\|_{L^2(\bbR^2)} \|A\|^2_{H^2(\bbR^2)} \\[2mm]
& = & C \|A\|_{L^2(\bbR^2)} \|A\|^3_{H^2(\bbR^2)}
\end{eqnarray*}
(we have used the fact that $H^2(\bbR^2)$ is an algebra and the embedding 
$H^1(\bbR^2)\subset L^q(\bbR^2)$ for all $q>2$). 

Similarly, for any $p>2$
\begin{equation*}
\left\|\frac{\partial_x(|A|^2)}{1+|A|^2}\right\|_{L^p(\bbR^2)}
\leq 2 \|A\|_{L^{2p}(\bbR^2)} \|\partial_x A\|_{L^{2p}(\bbR^2)} \\[2mm]
\leq C \|A\|_{H^1(\bbR^2)} \|A\|_{H^2(\bbR^2)}.
\end{equation*}
Finally for any $p>2$ 
\begin{equation*}
\|F\|_{L^p(\bbR^2)} 
\leq C \|A\|^2_{H^2(\bbR^2)} (1+\|A\|_{L^2(\bbR^2)}\|A\|_{H^2(\bbR^2)}).  
\end{equation*}
and by elliptic regularity
\begin{equation*}
\|\nabla\ph\|_{W^{1,p}(\bbR^2)} 
\leq C \|A\|^2_{H^2(\bbR^2)} (1+\|A\|_{L^2(\bbR^2)}\|A\|_{H^2(\bbR^2)}),
\hspace{1cm} \forall p>2.   
\end{equation*}
We now check that 
$\nabla\ph\cdot\nabla|A|^2/(1+|A|^2)\in H^1(\bbR^2)$. This easily reduces to
showing that $\nabla(\nabla\ph\cdot\nabla|A|^2)\in L^2(\bbR^2)$. For 
$(i,j)\in\{1,2\}$, $\partial_{x_i}\partial_{x_j}\ph\in L^p(\bbR^2)$ since 
\begin{equation*}
\widehat{\partial_{x_i}\partial_{x_j}\ph} =
\frac{\xi_i\xi_j}{|\xi|^2}\widehat{\Delta\ph} 
\textrm{ and } \Delta\ph\in L^p(\bbR^2),\ p>2. 
\end{equation*}
Thus $\partial_{x_i}\partial_{x_j}\ph \nabla|A|^2 \in L^2(\bbR)$
($\nabla|A|^2\in H^1(\bbR^2)\subset L^q(\bbR^2)$, $\forall q>2$). 

On the other hand, taking $p>2$ we see that $\nabla\ph\in L^\infty(\bbR^2)$ and
thus $\nabla\ph \partial_{x_i}\partial_{x_j} |A|^2 \in L^2(\bbR^2)$. 

It is also easy to check that $\partial_x(|A|^2)/(1+|A|^2)\in H^1(\bbR^2)$. 

\noindent
Finally, $\Delta\ph = F \in H^1(\bbR^2)$, proving that 
$\nabla\ph\in H^2(\bbR^2)$ with an estimate of the form 
\begin{equation*}
\|\nabla \ph\|_{H^2(\bbR^2)} \leq P(\|A\|_{H^2(\bbR^2)}),
\end{equation*} 
where $P$ is a polynomial vanishing at 0, which proves \eqref{Eq_P}.
\endproof

\remark
All above estimates are therefore uniform in time, and if
$A\in\calC([0,T];H^2(\bbR^2))$ for some $T>0$, one has
\begin{equation*}
\|\nabla \ph\|_{\calC([0,T];H^2(\bbR^2))} 
\leq P(\|A\|_{\calC([0,T];H^2(\bbR^2))}).
\end{equation*}

\subsection{Solitary waves --- non existence results}

We now look for solitary wave solutions of \eqref{Eq_Adiv}, that is solutions of 
the form $(e^{i\omega t}U(x),\phi(x))$ with $x\in\bbR^d$, $\omega\in\bbR$, 
$U\in H^1(\bbR^d)$, and $\phi\in H$. Thus $(U,\phi)$ should satisfy the
system
\begin{equation}
\label{Eq_ZASW}
\left\{
\begin{array}{l}
- \Delta U + \omega U = a U \partial_x\phi, \\[2mm]
\div ( (1+|U|^2) \nabla \phi ) = \partial_x(|U|^2).
\end{array}
\right.
\end{equation}
The existence of nontrivial solutions of \eqref{Eq_ZASW} is an open problem. 
Note that \eqref{Eq_ZASW} does not seem to be the Euler--Lagrange equation 
associated to a variational problem. We have however:
\begin{proposition}
{\em (i)} Let $a=-1$ {\em (}defocusing case{\em )}. Then no nontrivial solution  
of \eqref{Eq_ZASW} exists for $\omega\geq0$. \\
{\em (ii)} Let $a=1$ {\em (}focusing case{\em )}. No nontrivial solution of 
\eqref{Eq_ZASW} exists for $\omega\geq1$. \\
{\em (iii)} Let $a=\pm1$. No nontrivial solution of \eqref{Eq_ZASW} exists 
if $\omega<0$ provided $\partial_x\phi = O(1/|\bfx|^{1+\eps})$, $\eps>0$ as 
$|\bfx|\to+\infty$.
\end{proposition}

\proof
From \eqref{Eq_ZASW} we have
\begin{equation*}
\int_{\bbR^d} |\nabla U|^2 d\bfx + \omega \int_{\bbR^d} |U|^2 d\bfx 
= a \int_{\bbR^d} |U|^2 \partial_x\phi d\bfx,
\end{equation*}
\begin{equation*}
\int_{\bbR^d} (1+|U|^2)|\nabla\phi|^2 d\bfx 
= - \int_{\bbR^d} |U|^2 \partial_x\phi d\bfx,
\end{equation*}
and 
\begin{equation}
\label{Eq_ZASWH1Claw}
\int_{\bbR^d} |\nabla U|^2 d\bfx + \omega \int_{\bbR^d} |U|^2 d\bfx 
- a \int_{\bbR^d} (1+|U|^2)|\nabla\phi|^2 d\bfx = 0,
\end{equation}
which proves (i). Now independent of the sign of $a$, and from 
\eqref{Eq_ZAH1Est} and \eqref{Eq_ZASWH1Claw},
\begin{equation*}
\int_{\bbR^d} |\nabla U|^2 d\bfx + \omega \int_{\bbR^d} |U|^2 d\bfx 
\leq \int_{\bbR^d} |U|^2 d\bfx. 
\end{equation*}
Thus
\begin{equation*}
\int_{\bbR^d} |\nabla U|^2 d\bfx + (\omega-1) \int_{\bbR^d} |U|^2 d\bfx \leq0,
\end{equation*}
which proves (ii). Part (iii) results from \refcite{kato}.
\endproof

\subsection{The Cauchy problem}

We consider the system 
\begin{equation}
\label{Eq_ZACauchy}
\left\{
\begin{array}{l}
i\partial_t A + \Delta A = -aA\partial_x\ph, \\[2mm]
\div \left((1+|A|^2) \nabla \ph\right) = \partial_x(|A|^2), \\[2mm]
A(\cdot,0) = A_0.
\end{array}
\right.
\end{equation}

\begin{theorem}
\label{Th_ZACauchy}
Let $A_0\in H^2(\bbR^2)$. \\
Then there exists $T_0>0$ and a unique solution $(A,\nabla\ph)$ of 
\eqref{Eq_ZACauchy} such that $A\in\calC([0,T_0];H^2(\bbR^2))$ and 
$\nabla\ph\in\calC([0,T_0];H^2(\bbR^2))$. Moreover
\begin{equation*}
\|A(\cdot,t)\|_{L^2(\bbR^2)} = \|A_0\|_{L^2(\bbR^2)}, \hspace{1cm}0\leq t\leq T_0
\end{equation*}
and
\begin{equation*}
\int_{\bbR^2} (1+\frac12|A|^2) |\nabla\ph|^2 d\bfx 
\leq \frac12 \int_{\bbR^2} |A_0|^2 d\bfx, \hspace{1cm}0\leq t\leq T_0.
\end{equation*}

\end{theorem}

\proof
\textit{Uniqueness.} Let $(A,\nabla\ph)\in L^\infty(0,T;H^2(\bbR^2))$ and  
$(B,\nabla\psi)\in L^\infty(0,T;H^2(\bbR^2))$ two solutions of
\eqref{Eq_ZACauchy} with $A(\cdot,0)=B(\cdot,0)$. Then from
\eqref{Eq_ZACauchy}$_2$ one gets
\begin{equation*}
\Delta (\ph-\psi) + \div(|A|^2\nabla\ph-|B|^2\nabla\psi)
= \partial_x(|A|^2) - \partial_x(|B|^2),
\end{equation*}
yielding
\begin{equation}
\label{Eq_ZACauchy_unique}
\begin{aligned}
\int_{\bbR^2} & |\nabla (\ph-\psi)|^2 d\bfx 
+ \int_{\bbR^2} |A|^2 |\nabla(\ph-\psi)|^2  d\bfx \\
& = \int_{\bbR^2} (|A|^2-|B|^2)\partial_x (\ph-\psi) d\bfx 
- \int_{\bbR^2} (|A|^2-|B|^2) \nabla\psi\cdot\nabla(\ph-\psi) d\bfx .
\end{aligned}
\end{equation}
Observing that $|A|^2-|B|^2=A(A-\bar B)+\bar B(A-B)$, the R.H.S. of
\eqref{Eq_ZACauchy_unique} is majorized by
\begin{equation*}
\begin{aligned}
\frac14 \int_{\bbR^2} |\partial_x & (\ph-\psi)|^2 d\bfx 
+ (\|A\|_{L^\infty(\bbR^2)}+\|B\|_{L^\infty(\bbR^2)}) 
\int_{\bbR^2} |A-B|^2 d\bfx \\
& + \frac14 \int_{\bbR^2} |\nabla (\ph-\psi)|^2 d\bfx \\
& + (\|A\|_{L^\infty(\bbR^2)}+\|B\|_{L^\infty(\bbR^2)})\|
\nabla\psi\|_{L^\infty(\bbR^2)} \int_{\bbR^2} |A-B|^2 d\bfx 
\end{aligned}
\end{equation*}

\noindent
and by Sobolev embedding
\begin{equation}
\label{Eq_ZACauchy_nablapsi}
\|\nabla(\ph-\psi)\|_{L^2(\bbR^2)} 
\leq C(\|A\|_{H^2(\bbR^2)},\|B\|_{H^2(\bbR^2)},\|\nabla\psi\|_{H^2(\bbR^2)})
\|A-B\|_{L^2(\bbR^2)}.
\end{equation}
On the other hand, we obtain readily from \eqref{Eq_ZACauchy}$_1$ that
\begin{equation*}
\frac12 \frac{d}{dt} \int_{\bbR^2} |A-B|^2 d\bfx 
\leq \int_{\bbR^2} |A-B|^2 |\partial_x\ph| d\bfx
+ \int_{\bbR^2} |B||\partial_x(\ph-\psi)||A-B| d\bfx
\end{equation*}
which together with \eqref{Eq_ZACauchy_nablapsi} and the Cauchy-Schwarz lemma
yields
\begin{equation*}
\begin{aligned}
\frac12 \frac{d}{dt} & \int_{\bbR^2} |A-B|^2 d\bfx \\
& \leq C(\|A\|_{H^2(\bbR^2)},\|B\|_{H^2(\bbR^2)},\|\nabla\ph\|_{H^2(\bbR^2)},
\|\nabla\psi\|_{H^2(\bbR^2)}) \|A-B\|_{L^2(\bbR^2)}
\end{aligned}
\end{equation*}
and $A=B$ by Gronwall lemma. \\[2mm]
\textit{$H^2$ a priori estimate.} We derive a (formal) $H^2$ 
\textit{a priori} estimate on the solution of \eqref{Eq_ZACauchy}. Since 
$H^2(\bbR^2)$ is an algebra, we deduce from Lemma~\ref{Lm_Invdivgrad} that 
\begin{equation}
\label{Eq_EstH2_Aphx}
\|A\partial_x \ph\|_{\calC([0,T];H^2(\bbR^2))} 
\leq \|A\|_{\calC([0,T];H^2(\bbR^2))} P(\|A\|_{\calC([0,T];H^2(\bbR^2))}), 
\end{equation}
where $P$ was introduced in \eqref{Eq_P}. From the energy estimate
\begin{equation*}
\frac12 \frac{d}{dt} \|A(\cdot,t)\|^2_{H^2(\bbR^2)} 
\leq C \|A\partial_x \ph (\cdot,t)\|_{H^2(\bbR^2)}\|A(\cdot,t)\|_{H^2(\bbR^2)},
\end{equation*}
we infer with \eqref{Eq_EstH2_Aphx} the local $H^2$ bound
\begin{equation}
\label{Eq_EstH2_A}
\|A(\cdot,t)\|_{H^2(\bbR^2)} \leq C \left(\|A_0\|_{H^2(\bbR^2)} \right)
\textrm{ for } 0<t<T_0, 
\end{equation}
$T_0<T$ sufficiently small. \\[2mm]
\textit{Approximation of \eqref{Eq_ZACauchy}.} The strategy is now to implement
a compactness method using the (justified) \textit{a priori} estimate 
\eqref{Eq_EstH2_A}. For $\eps>0$, we consider the system
\begin{eqnarray}
\label{Eq_ZACauchy_epsA}
i\partial_t A^\eps + \Delta A^\eps & = & -aA^\eps\partial_x\ph^\eps, \\[2mm]
\label{Eq_ZACauchy_epsph}
\div \left((1+\eps \Delta^2+ |A^\eps|^2) \nabla \ph^\eps\right) 
& = & \partial_x(|A^\eps|^2), \\[2mm]
\label{Eq_ZACauchy_epsA0}
A^\eps(\cdot,0) & = & A_0. 
\end{eqnarray}
Solving $\nabla\ph^\eps$ in terms of $A^\eps$, we obtain from
\eqref{Eq_ZACauchy_epsph} that $\nabla\ph^\eps$ satisfies
\begin{equation}
\label{Eq_ZACauchy_uniform}
\eps \int_{\bbR^2} |\Delta\nabla\ph^\eps|^2 d\bfx
+ \int_{\bbR^2} (1+\frac12|A^\eps|^2)|\nabla\ph^\eps|^2 d\bfx
\leq \frac12 \int_{\bbR^2} |A^\eps|^2 d\bfx.
\end{equation}

\noindent
\textit{Well-posedness of approximate system.} We now check that the Cauchy 
problem \eqref{Eq_ZACauchy_epsA}--\eqref{Eq_ZACauchy_epsA0} is globally 
well-posed in $H^2(\bbR^2)$. Let first $A^\eps$, $B^\eps\in H^2(\bbR^2)$ and 
$\ph^\eps$, $\psi^\eps$ the corresponding solutions of
\eqref{Eq_ZACauchy_epsph}. Proceeding as in the uniqueness proof above, one 
gets 
\begin{equation}
\label{Eq_ZACauchy_unifeps}
\begin{aligned}
\eps \int_{\bbR^2} & |\nabla\Delta (\ph^\eps-\psi^\eps)|^2 d\bfx 
+ \int_{\bbR^2} |\nabla(\ph^\eps-\psi^\eps)|^2 d\bfx \\
& \leq C(\|A^\eps\|_{H^2(\bbR^2)},\|B^\eps\|_{H^2(\bbR^2)},
\|\nabla\ph^\eps\|_{H^2(\bbR^2)},\|\nabla\psi^\eps\|_{H^2(\bbR^2)}) 
\|A^\eps-B^\eps\|_{L^2(\bbR^2)}.
\end{aligned}
\end{equation}
Denoting $\partial_x\ph^\eps$ by $F^\eps(A^\eps)$ we write 
\eqref{Eq_ZACauchy_epsA} on the Duhamel form with $S(t)=\exp(it\Delta)$,
\begin{equation}
\label{Eq_ZACauchy_Duhameleps}
A^\eps(t) = S(t) A_0 - a \int_0^t S(t-s)A^\eps F^\eps(A^\eps) ds.
\end{equation}
Using \eqref{Eq_ZACauchy_unifeps} and the unitarity of $S(t)$ in $H^s(\bbR^2)$,
we deduce that the R.H.S. of \eqref{Eq_ZACauchy_Duhameleps} defines a
contraction in $\calC([0,T_\eps];H^2(\bbR^2))$ for some $T_\eps>0$. 

This implies the local well-posedness of 
\eqref{Eq_ZACauchy_epsA}--\eqref{Eq_ZACauchy_epsA0} in $H^2(\bbR^2)$. Using the 
$H^2$ bound \eqref{Eq_ZACauchy_unifeps} on $\nabla\ph^\eps$, we infer from
\eqref{Eq_ZACauchy_epsA} an \textit{a priori} bound in 
$\calC([0,T];H^2(\bbR^2))$ for $A^\eps$ and for all $T>0$. This proves that the
Cauchy problem \eqref{Eq_ZACauchy_epsA}--\eqref{Eq_ZACauchy_epsA0} is globally
well-posed, for any fixed $\eps>0$.\\[2mm]
\textit{Limit $\eps\to0$.} Now we have the bounds \eqref{Eq_ZACauchy_uniform} 
and 
\begin{equation}
\label{Eq_ZACauchy_boundeps}
\|\nabla\ph^\eps(\cdot,t)\|_{L^2(\bbR^2)} 
+ \sqrt\eps \|\Delta\nabla\ph^\eps(\cdot,t)\|_{L^2(\bbR^2)} \leq C, 
\hspace{1cm} 0\leq t\leq T,
\end{equation}
where $C$ and $T$ do not depend on $\eps$. Moreover, from
\eqref{Eq_ZACauchy_epsA} and \eqref{Eq_ZACauchy_boundeps} we have a bound on
$\partial_t A^\eps$ which is independent of $\eps$: 
\begin{equation*}
\|\partial_t A^\eps(\cdot,t)\|_{L^2(\bbR^2)} \leq C, 
\hspace{1cm} 0\leq t\leq T.
\end{equation*}
It is now standard to pass to the limit as $\eps\to0$ (see \refcite{lions}). 
By the Aubin--Lions compactness lemma, we obtain a subsequence
$(A^\eps,\nabla\ph^\eps)$ such that $A^\eps\to A$ in $L^\infty(0,T;H^2(\bbR^2))$
weak-star and $L^2(0,T;H^1_{\textrm{loc}}(\bbR^2))$ strongly, 
$\nabla\ph^\eps\to \nabla\ph$ in $L^\infty(0,T;H^2(\bbR^2))$
weak-star and $L^2([0,T]\times\bbR^2)$ weakly. The limit 
$(A,\nabla\ph)$ belongs to $(L^\infty(0,T;H^2(\bbR^2)))^2$ and satisfies
\eqref{Eq_ZACauchy}. In fact \eqref{Eq_ZACauchy}$_1$ is satisfied in
$L^2(\bbR^2)$ and \eqref{Eq_ZACauchy}$_2$ is satisfied in $H^1(\bbR^2)$.

The fact that $(A,\nabla\ph)\in (\calC(0,T;H^2(\bbR^2)))^2$ results from the
Bona--Smith approximation (see \refcite{bona-smith}).
\endproof

\remark 
We do not know whether the local solution obtained in Theorem \ref{Th_ZACauchy}
is global or not.

\section{Conclusion}

We have given a full description of how to derive from the Kukhtarev equations
an asymptotic model for the propagation of light in a photorefractive
medium. This derivation is only heuristic insofar as asymptotics are not
justified, which would be out of reach now. Some properties of photorefractive
media such as memory have also been neglected.

The 1D asymptotic model is a saturated nonlinear Schr\"odinger equation the 
Cauchy problem of which is studied (in any space dimension) in $L^2$ and $H^1$. 
We also prove the existence of solitary waves in one and higher dimensions. An 
interesting and open issue would be to study the transverse stability of the 
1D solitary waves in the framework of the asymptotic model.

For the 2D asymptotic model (the Zozulya--Anderson model) we also have studied 
the Cauchy problem and the non-existence of solitary waves. The question of
imposing other boundary conditions, not vanishing in one space direction, can 
also be addressed to treat a wider range of experimental applications.   

\appendix

\section{Non-Existence of Solitary Waves in Non-Physical Cases}

The goal is here to complete the results of Corollary \ref{Cor_NLSnoSW} for 
$\omega<0$ with no decaying assumption. We have already seen that Eq.
\eqref{Eq_NLSSWPohozaev} implies that no solitary wave may exist for $d=1,2$ and
$a=1$ (focusing case). 

To go further, let us use both Eqs. \eqref{Eq_NLSSWEnergy} and
\eqref{Eq_NLSSWPohozaev} to obtain
\begin{equation*}
\int_{\bbR^d} \left(2\omega + \frac{(d-2)a|U|^2}{1+|U|^2}\right) |U|^2 d\bfx 
- ad  \int_{\bbR^d} \left[ |U|^2 - \ln (1+|U|^2) \right] d\bfx = 0.
\end{equation*}
We set 
\begin{equation*}
F(X) = \left(2\omega + \frac{(d-2)aX}{1+X}\right) X - ad (X - \ln (1+X)),
\end{equation*}
and we know that $\int_{\bbR^d} F(|U|^2) d\bfx = 0$. Now $F(0)=0$ and 
\begin{equation*}
F'(X)= \frac{2X^2(\omega-a)+X(4\omega-(4-d)a)+2\omega}{(1+X)^2} < 0, 
\end{equation*}
if $\omega<0$, $a=1$ and $d=3,4$.
Therefore $F(|U|^2)=0$ a.e. By a bootstrapping argument, we notice that any
$H^1$ solution to Eq. \eqref{Eq_NLSSW} is indeed in $H^k$ for all $k$ and
therefore continuous. Hence $F(|U|^2)=0$ on $\bbR^d$. Since $F'(X)<0$ the only
possible value for $U$ is $U=0$ on $\bbR^d$.

We can refine this result, finding other parameter ranges for which 
$2X^2(\omega-a)+X(4\omega-(4-d)a)+2\omega<0$. If $d\geq5$ and $a=1$, this holds
for $\omega\leq -(d-4)/4$. Moreover,
\begin{equation*}
\begin{aligned}
2X^2(\omega-a)&+X(4\omega-(4-d)a)+2\omega \\
& = 2\left[(\omega-a)(X-1)^2+X(4\omega-4a+\frac d2 a) + \omega+a\right].
\end{aligned}
\end{equation*}
No solitary wave can exist for $a=-1$ and $\omega\leq-1$. Hence we complete 
Corollary \ref{Cor_NLSnoSW} with

\def\thecorollary{A.\arabic{corollary}}
\setlength\leftmargini{2pc}
\begin{corollary}
No non-trivial solitary wave (solution of \eqref{Eq_NLSSW}) of the saturated NLS
equation exists when
\begin{itemize}
\item[{\em (i)}] $a=-1$ {\em (}defocusing case{\em )}, for $\omega\leq-1$.
\item[{\em (ii)}] $a=1$ {\em (}focusing case{\em )}, for $\omega\leq0$, if $d=3,4$ and
  $\omega\leq-(d-4)/4$ if $d\geq5$.
\end{itemize}
\end{corollary}

% \section*{References}


{\small
\begin{thebibliography}{10}

\bibitem{abramowitz-steglun}
M. Abramowitz and I. A. Steglun,
\newblock {\em Handbook of Mathematical Functions}
(National Bureau of Standards, 1964).

\bibitem{berestycki-lions-peletier}
H.~Berestycki, P.-L. Lions and L. A. Peletier,
\newblock An {ODE} approach to the existence of positive solutions for
  semilinear problems in {$\bbR^N$},
\newblock {\em Indiana Univ. Math. J.} \textbf{30} (1981) 141--157.

\bibitem{bona-smith}
J. L. Bona and R. Smith,
\newblock The initial value problem for the {K}orteweg--de {V}ries equation,
\newblock {\em Philos. Trans. Royal. Soc. London A} \textbf{278} (1975) 
555--601.

\bibitem{cazenave1}
T. Cazenave,
\newblock {\em An Introduction to Nonlinear {S}chr\"odinger Equations,} 
3rd edition, Textos de Metodos Matematicos. Vol. 26 (Instituto de Matematica, 
Universidade Federal do Rio de Janeiro, 1996).

\bibitem{cazenave2}
T. Cazenave,
\newblock {\em Semilinear {S}chr\"odinger Equations}, 
Courant Lecture Notes in Mathematics. Vol. 10 (Courant Institute of 
Mathematical Sciences, New York, 2003).

\bibitem{delre-crosignani-diporto}
E. DelRe, B. Crosignani and P. {Di Porto},
\newblock Photorefractive spatial solitons,
\newblock in {\em Spatial Solitons,} S. Trillo and W. Torruellas (eds), 
Springer Series in Optical Sciences (Springer, 2001), pp. 61--85.

\bibitem{ghidaglia-saut1}
J.-M. Ghidaglia and J.-C. Saut,
\newblock On the initial value problem for the {D}avey--{S}tewartson systems,
\newblock {\em Nonlinearity} \textbf{3} (1990) 475--506.

\bibitem{kato}
T. Kato,
\newblock Growth properties of solutions of the reduced wave equation with a
  variable coefficient,
\newblock {\em Commun. Pure Appl. Math.} \textbf{12} (1959) 403--425.

\bibitem{kukhtarev-markow-odoluv-soskin-vinetskii}
N. V. Kukhtarev, V. B. Markow, S. G. Odoluv, M. S. Soskin  and V. L. Vinetskii,
\newblock Holographic storage in electrooptic crystals,
\newblock {\em Ferroelectrics} \textbf{22} (1979) 949--960.

\bibitem{lions}
J.-L. Lions,
\newblock {\em Quelques m\'ethodes de r\'esolution des probl\`emes aux limites
  non lin\'eaires} (Dunod, Paris, 1969).

\bibitem{mamaev-saffman-zozulya2}
A. V. Mamaev, M. Saffman and A. A. Zozulya,
\newblock Break-up of two-dimensional bright spatial solitons due to transverse
  modulational instability,
\newblock {\em Europhys. Lett.} \textbf{35} (1996) 25--30.

\bibitem{mamaev-saffman-zozulya1}
A. V. Mamaev, M. Saffman  and A. A. Zozulya,
\newblock Propagation of dark stripe beams in nonlinear media: {S}nake
  instability and creation of optical vertices, 
\newblock {\em Phys. Rev. Lett.} \textbf{76} (1996) 2262--2265.

\bibitem{stepken-kaiser-belic-krolikowski}
A. Stepken, F. Kaiser, M. R. Beli\'c  and W. Kr\'olikowski,
\newblock Interaction of incoherent two-dimensional photorefractive solitons,
\newblock {\em Phys. Rev. E} \textbf{58} (1998) R4112--R4115.

\bibitem{tikhonenko-christou-luther-davies}
V. Tikhonenko, J. Christou and B. Luther-Davies,
\newblock Three-dimensional bright spatial soliton collision and fusion in a
  saturable nonlinear medium,
\newblock {\em Phys. Rev. Lett.} \textbf{76} (1996) 2698--2701.

\bibitem{wolfersberger-fressengeas-maufroy-kugel}
D. Wolfersberger, N. Fressengeas, J. Maufroy and G. Kugel,
\newblock Self-focusing of a single laser pulse in a photorefractive medium,
\newblock {\em Phys. Rev. E} \textbf{62} (2000) 8700--8704.

\bibitem{zozulya-anderson}
A. A. Zozulya and D. Z. Anderson,
\newblock Propagation of an optical beam in a photorefractive medium in the
  presence of a photogalvanic nonlinearity or an externally applied electric
  field,
\newblock {\em Phys. Rev. A} \textbf{51} (1995) 1520--1531.

\end{thebibliography}
}
\end{document}